\documentstyle{amsppt}
\magnification=\magstep1
\hcorrection{0in}
\vcorrection{0in}
\pagewidth{6.6truein}
\pageheight{9truein}
\NoRunningHeads
\def\Q{\Bbb Q}
\def\C{\Bbb C}
\def\O{\Cal O}
\def\rup#1{\left\lceil #1 \right\rceil}
\def\rdown#1{\left\lfloor #1 \right\rfloor}
\def\proper#1{\widetilde{#1}}
\def\exc#1{{#1}^{\text{exc}}}
\def\Supp{\mathop{Supp}}
\topmatter
\title Effective base point freeness on a normal surface
\endtitle
\author Takeshi Kawachi \endauthor
\affil Department of Mathematics, Tokyo Institute of Technology \endaffil
\address 2-12-1 Oh-okayama, Meguro-ku, Tokyo, JAPAN, 152 \endaddress
\email kawachi\@math.titech.ac.jp \endemail
\abstract
We prove effective base point freeness of the adjoint linear system
on normal surfaces with a boundary.
\endabstract
\endtopmatter
\document
\head 0. Introduction \endhead
Let $Y$ be a compact normal two dimensional projective variety over $\C$
(we call it ``normal surface'' for short).
Let $y\in Y$ be a given point, and $D$ be a nef and big $\Q$-divisor
on $Y$ such that $K_Y+\rup{D}$ is a Cartier divisor, 
There is various numerical criteria on $D$ when the adjoint linear system
$|K_Y+\rup{D}|$ is free at $y$.
The main theorem includes all earlier result and gives the first effective
version if $y$ is log-terminal singularity.
\par
To describe the main theorem we prepare notation.
\par
\definition{Definition}
Let $Y$ be a normal surface, and $B$ be an effective $\Q$-divisor on $Y$.
Let $y$ be a point on $Y$.
Then the triple $(Y,B,y)$ is said to be a germ of {\it quasi-log-terminal}
singularity if the following conditions are satisfied:
\roster
\item $[B]=0$, where $[B]$ is the integral part of $B$.
\item Let $f\:X\to Y$ be the minimal resolution of $y$
      (if $y$ is a smooth point then let $f$ be the blowing-up at $y$).
      If we write
      $$
         K_X=f^*(K_Y+B)+\sum a_jF_j\quad\text{(for $a_j\in\Q$)}
      $$
      then $a_j>-1$ whenever $F_j$ is $f$-exceptional.
\endroster
\enddefinition
\par
\remark{Remark}
The differences with ``log-terminal'' are:
\roster
\item We don't mention where $K_Y+B$ is $\Q$-Cartier or not.
      The pull-back and the intersection are used by Mumford's $\Q$-valued
      pull-back and intersection theory.
\item $f$ is only ``minimal resolution''.
      Thus the $\Q$-divisor $\sum a_jF_j$ may not be normal crossings.
      Hence ``quasi-log-terminal'' may not be ``log-terminal''
      (cf\. [KMM]).
\item If $(Y,B,y)$ is quasi-log-terminal then $(Y,0,y)$ which is ignored
      the boundary, is log-terminal.
\endroster
\endremark
\par
Let $(Y,B,y)$ be a germ of normal surface singularity with a boundary $B$,
and let $f\:X\to Y$ be the minimal resolution of $y$
(if $y$ is a smooth point then the blowing-up at $y$).
Let $Z$ is the fundamental cycle of $y$
(if $y$ is a smooth point then the exceptional divisor of $f$).
Let $\Delta_y=f^*K_Y-K_X$.
$\Delta_y$ is an effective $\Q$-divisor if $y$ is singular,
and $\Delta_y=-Z$ if $y$ is smooth.
\par
\definition{Definition}
$$
  \delta_y=
    \cases
      -(Z-\Delta_y)^2, & \text{if $y$ is quasi-log-terminal} \\
      0, & \text{if $y$ is not quasi-log-teminal}
    \endcases
$$
\enddefinition
\par
\definition{Definition}
We define
$$
  \mu(B,y)=\max\{\mu\mid \mu(Z-\Delta_y)\leq f^*B\}
$$
\enddefinition
\par
\remark{Remark}
\roster
\item This $\delta_y$ satisfies $0\leq\delta_y\leq 4$ ([KM]).
\item If $y$ is a smooth point then $2\mu=\text{mult}_y B$.
      If $y$ is quasi-log-terminal then $\mu<1$.
\endroster
\endremark
\par
The main result is the following.
\proclaim{Theorem}
Let $(Y,B)$ be a normal surface with a boundary $B$ and assume that $[B]=0$.
Let $M$ be a Cartier divisor on $Y$ such that $M-(K_Y+B)$ is nef big.
Let $D$ be an effective $\Q$-divisor on $Y$ which is $\Q$-linearly equivalent
to $M-(K_Y+B)$.
Let $y$ be a point on $Y$.
Let $\mu=\mu(B,y)$.
Assume that $D^2>(1-\mu)^2\delta_y$.
If $DC\geq (1-\mu)\delta_y/2$ for all curves $C$ in $Y$ such that $y\in C$
then the linear system $|M|$ is free at $y$.
Furthermore if $y$ is a singular point which is not of type $A_n$ then
$|M|$ is free at $y$.
\endproclaim
\par
We can get the following immediately.
\proclaim{Corollary}
Let $Y$ be a normal surface.
Let $D$ be a nef $\Q$-divisor on $Y$ such that $K_Y+\rup{D}$ is Cartier.
Let $y$ be a point on $Y$.
Let $B=\rup{D}-D$ and $\mu=\mu(B,y)$
Assume that $D^2>(1-\mu)^2\delta_y$.
If $DC\geq (1-\mu)\delta_y/2$ for all curves $C$ in $Y$ such that $y\in C$
then the linear system $|K_Y+\rup{D}$ is free at $y$.
\endproclaim
\par
\remark{Remark}
If $y$ is a smooth point or a rational double point then
this theorem is essentially the theorem of Ein-Lazarsfeld ([EL]).
If $y$ is not a quasi-log-teminal point then the result are proved in [ELM].
If $B=0$ then the result is described in [KM].
\endremark
In section 1, we remark the facts about quasi-log-terminal singularity.
These are used by the proof of the main theorem.
In section 2, we prove the theorem.
\head 1. Facts about quasi-log-terminal singularity \endhead
Throughout this section,
let $(Y,B,y)$ be a germ of quasi-log-terminal surface singularity.
Let $\Delta_y=f^*K_Y-K_X=\sum a_j\Delta_j$
where $\Delta_i$ is an irreducible exceptional divisor lying over $y$.
Let $B=B_y+B'$ and
$f^*B_y=\proper{B_y}+\exc{B_y}=\proper{B_y}+\sum b_j\Delta_j$
where $B_y$ consists of the components which passing through $y$,
$B'$ is disjoint from $y$ and $\proper{B_y}$ is the proper
transformation of $B_y$.
Let $Z$ be the fundamental cycle of $y$.
Hence we have
$$
  f^*(K_Y+B)=K_X+\proper{B_y}+\Delta_y+\exc{B_y}+f^*B'.
$$
First we bound the value of $\mu(B,y)$.
we recall the definition of $\mu(B,y)$.
$$
  \mu(B,y)=\max\{\mu\mid \mu(Z-\Delta_y)\leq f^*B\}
$$
\par
\proclaim{Lemma 1}
If $(Y,B,y)$ is quasi-log-terminal then $\mu(B,y)<1$.
\endproclaim
\par
\demo{Proof}
If $\mu(B,y)\geq 1$ then $\exc{B_y}\geq Z-\Delta_y$.
Hence we have $\Delta_y+\exc{B_y}\geq Z$.
On the other hand, since $(Y,B,y)$ is quasi-log-terminal,
all coefficients of $\Delta_y+\exc{B_y}$ must be less than~$1$.
That is contradiction. \qed
\enddemo
\par
Since $B$ is an effective divisor, $\exc{B_y}$ is also effective.
Therefore all coefficients of $\Delta_y$ are less than $1$.
This is the case of classical ``log-terminal'' surface singularity
without boundary.
\par
We defined $\delta_y=-(Z-\Delta_y)^2$.
This is bounded by the following lemma (cf\. [KM]).
\par
\proclaim{Lemma 2}
$0<\delta_y<2$ if $y$ is log-terminal but not a smooth or
a rational double point.
\endproclaim
\par
\demo{Proof}
Since the intersection matrix $\|\Delta_i\Delta_j\|$ is negative definite
and $Z-\Delta_y$ is effective and non-zero, we have $\delta_y>0$.
On the other hand, we have
$$\align
  (Z-\Delta_y)^2 &=Z(Z-\Delta_y)-\Delta_y(Z-\Delta_y) \\
  &= Z(Z+K_X)+K_X(Z-\Delta_y) \\
  &> Z(Z+K_X) = 2\mathop{Pa}(Z)-2.
\endalign$$
Since the log-terminal singularity is the rational singularity,
we have $\mathop{Pa}(Z)=0$ by [A].
Therefore we have $\delta_y<2$. \qed
\enddemo
\par
It is easy to prove that $\delta_y=4$ if $y$ is a smooth point and
that $\delta_y=2$ if $y$ is a rational double point.
\par
Let $C$ be an irreducible divisor on $Y$ and $y\in C$.
Let $f^*C=\proper{C}+\sum c_{i,j}\Delta_j$, and
$\proper{C}\Delta_i=1$ for an $i$ and $\proper{C}\Delta_j=0$ for $j\ne i$.
That is, we assume that $\proper{C}$ meets only $\Delta_i$.
Next, we study the value $a_i+c_{i,i}$.
Since all log-terminal singularities are classified in [B],
we examine it in each cases.
we set indexes as in Appendix.
\par
\subhead 1.1 The case of type $A_n$ \endsubhead
Let $A=A(w_1,\dots,w_n)$ be the intersection matrix of type $A_n$
where $w_i=-\Delta_i^2$.
$$
  A(w_1,\dots,w_n)=
    \pmatrix
    -w_1 & 1      & 0      & 0    \\
    1    & -w_2   & \ddots & 0    \\
    0    & \ddots & \ddots & 1    \\
    0    & 0      &  1     & -w_n \\
    \endpmatrix.
$$
Let $a(w_1,\dots,w_n)=\det A(w_1,\dots,w_n)$, and
we define $a()=1$ for convention.
\par
Let $\widetilde{A}_{ij}$ be the $(i,j)$-component of $A^{-1}$.
Suppose $j<i$ then we have
$$\align
  \widetilde{A}_{ij}
  &= \frac1{\det A}(-1)^{i+j}
     \vmatrix
  -w_1 & 1      & 0        &   &          &          &      &        & 0    \\
  1    & \ddots & 1        & 0 &          &          &      &        &      \\
  0    & 1      & -w_{j-1} & 1 & 0        & 0        &      &        &      \\
       & 0      & 0        & 1 & -w_{j+1} & 1        & 0    &        &      \\
       &        &          & 0 & \ddots   & -w_{i-1} & 0    &        &      \\
       &        &          &   & 0        & 1        & 1    & 0      &      \\
       &        &          &   &          & 0        & -w_i & \ddots & 0    \\
       &        &          &   &          & 0        & 1    & \ddots & 1    \\
       & 0      &          &   &          &          & 0    & 1      & -w_n \\
     \endvmatrix \\
  &= \frac1{\det A}(-1)^{i+j}a(w_1,\dots,w_{j-1}) a(w_{i+1},\dots,w_n).
\endalign$$
If $j>i$ then we have the following as same as above.
$$
  \widetilde{A}_{ij}=
    \frac1{\det A}(-1)^{i+j}a(w_1,\dots,w_{i-1}) a(w_{j+1},\dots,w_n).
$$
Also we have
$$
  \widetilde{A}_{ii}=
    \frac1{\det A}a(w_1,\dots,w_{i-1}) a(w_{i+1},\dots,w_n).
$$
Since
$$
  A(w_1,\dots,w_n)
  \pmatrix a_1 \\ \vdots \\ a_n \endpmatrix
  = \pmatrix 2-w_1 \\ \vdots \\ 2-w_n \endpmatrix,
$$
we have the following.
$$\align
  a_i &= \sum_{j=1}^{i-1} \widetilde{A}_{ij}(2-w_j) +
         \widetilde{A}_{ii}(2-w_i) +
         \sum_{j=i+1}^n \widetilde{A}_{ij}(2-w_j) \\
      &= \frac1{\det A}\bigl\{a(w_{i+1},\dots,w_n) 
         \bigl((-1)^{i+1}(2-w_1)+\cdots+
         (-1)^{2i-1}a(w_1,\dots,w_{i-2})(2-w_{i-1})\bigr)+ \\
      &\hbox to .4in{\hfil}
         a(w_1,\dots,w_{i-1}) a(w_{i+1},\dots,w_n) (2-w_i)+ \\
      &\hbox to .4in{\hfil}
         a(w_1,\dots,w_{i-1})\bigl(
         (-1)^{2i+1}a(w_{i+2},\dots,w_n)(2-w_{i+1})+
         \cdots+(-1)^{i+n}(2-w_n)\bigr)\bigr\}.
\endalign$$
Since $a(w_1,\dots,w_j)=-w_j a(w_1,\dots,w_{j-1})-a(w_1,\dots,w_{j-2})$,
we have
$$\align
  &(-1)^{i+1}(2-w_1)+(-1)^{i+2}a(w_1)(2-w_2)+\cdots +
   (-1)^{2i-1}a(w_1,\dots,w_{i-2})(2-w_{i-1}) \\
  &= (-1)^{i+1}\bigl\{
     2-w_1+(-1)(2a(w_1)+a()+a(w_1,w_2))+\cdots+ \\
  &\hbox to .4in{\hfil}
     (-1)^{i-2}(2a(w_1,\dots,w_{i-2})+a(w_1,\dots,w_{i-3})+
                 a(w_1,\dots,w_{i-1}))\bigr\} \\
  &= (-1)^{i+1}\bigl(1+(-1)^{i-2}(a(w_1,\dots,w_{i-2})+
                     a(w_1,\dots,w_{i-1}))\bigr) \\
  &= (-1)^{i+1}-a(w_1,\dots,w_{i-2})-a(w_1,\dots,w_{i-1}).
\endalign$$
We also have the following by same calculation.
$$\align
  &(-1)^{2i+1}a(w_{i+2},\dots,w_n)(2-w_{i+1})+\cdots+(-1)^{i+n}(2-w_n) \\
  &= -a(w_{i+1},\dots,w_n)-a(w_{i+2},\dots,w_n)+(-1)^{i+n}.
\endalign$$
Hence we have
$$\align
  a_i &= \frac1{\det A}\bigl\{
         a(w_{i+1},\dots,w_n)
           ((-1)^{i+1}-a(w_1,\dots,w_{i-2})-a(w_1,\dots,w_{i-1}))+ \\
      &\hbox to .4in{\hfil}
         a(w_1,\dots,w_{i-1}) a(w_{i+1},\dots,w_n) (2-w_i)+ \\
      &\hbox to .4in{\hfil}
         a(w_1,\dots,w_{i-1})
           ((-1)^{i+n}-a(w_{i+1},\dots,w_n)-a(w_{i+2},\dots,w_n))
         \bigr\}\\
      &= \frac1{\det A}\bigl\{
          (-1)^{i+1}a(w_{i+1},\dots,w_n)+(-1)^{i+n}a(w_1,\dots,w_{i-1}) \\
      &\hbox to .4in{\hfil}
          -a(w_1,\dots,w_{i-2})a(w_{i+1},\dots,w_n)
          -a(w_1,\dots,w_{i-1})a(w_{i+1},\dots,w_n)w_i \\
      &\hbox to .4in{\hfil}
          -a(w_1,\dots,w_{i-1})a(w_{i+2},\dots,w_n)\bigr\} \\
      &= 1+\frac1{\det A}\{
         (-1)^{i+1}a(w_{i+1},\dots,w_n)+(-1)^{i+n}a(w_1,\dots,w_{i-1})\}.
\endalign$$
\par
Since $(A(w_1,\dots,w_n)(c_{i,j}))_j=-1$ if $j=i$ and $=0$ if $j\ne i$,
we have 
$$\align
  c_{i,j} &= -\widetilde{A}_{ji} \\
          &= \frac1{\det A}(-1)^{i+j+1}
               \cases
                 a(w_1,\dots,w_{j-1})a(w_{i+1},\dots,w_n), & \text{if } j\leq i \\
                 a(w_1,\dots,w_{i-1})a(w_{j+1},\dots,w_n), & \text{if } j>i. \\
               \endcases
\endalign$$
Therefore we have the following.
$$\align
  a_i+c_{i,i}
  &= 1+\frac1{|\det A|}(
     -|a(w_{i+1},\dots,w_n)|-|a(w_1,\dots,w_{i-1})| + \\
  &\hbox to .4in{\hfil}
     |a(w_1,\dots,w_{i-1})| |a(w_{i+1},\dots,w_n)|) \\
  &= 1+\frac1{|\det A|}\bigl(
     (|a(w_1,\dots,w_{i-1})|-1)(|a(w_{i+1},\dots,w_n)|-1)-1\bigr)
\endalign$$
Hence we get the following proposition.
\par
\proclaim{Proposition 1}
Let $C$ as above. Assume $y$ is of type $A_n$.
If $n\geq 3$ and $i\ne 1,n$ then $a_i+c_{i,i}\geq 1$.\newline
If $i=1$ or $n$ then $a_i+c_{i,i}=1-1/|\det A|$ (for $n\geq 1$).
\endproclaim
\par
In the case of type $A_n$ we have $\delta_y=2-a_1-a_n$.
Indeed, since $(\Delta_y-Z)\Delta_i=2-w_i-(-w_i+2)=0$ if $i\ne 1,n$
and $(\Delta_y-Z)\Delta_i=1$ if $i=1$ or $n$,
we get the equation $\delta_y=2-a_1-a_n$.
\subhead 1.2 The case of type $D_n$ \endsubhead
Let $D=D(w_1,\dots,w_{n-2})$ be the intersection matrix of type $D_n$
where $w_i=-\Delta_i^2$.
$$
  D(w_1,\dots,w_n)=
    \pmatrix
    -w_1 & 1      & 0      &          &    & 0  \\
    1    & -w_2   & \ddots & 0        &    &    \\
    0    & \ddots & \ddots & 1        & 0  & 0  \\
         & 0      & 1      & -w_{n-2} & 1  & 1  \\
         &        & 0      & 1        & -2 & 0  \\
    0    &        & 0      & 1        & 0  & -2 \\
    \endpmatrix.
$$
Let $d(w_1,\dots,w_{n-2})=\det D(w_1,\dots,w_{n-2})$, and
we define $d()=4$ for convention.
\par
Let $\widetilde{D}_{ij}$ be the $(i,j)$-component of $D^{-1}$.
As same as of type $A_n$, we have
$$
  \widetilde{D}_{ij}=\frac1{\det D}
  \cases
  (-1)^{i+j}a(w_1,\dots,w_{j-1})d(w_{i+1},\dots,w_{n-2}),
    & i\leq n-2, j\leq i \\
  (-1)^{i+j}a(w_1,\dots,w_{i-1})d(w_{j+1},\dots,w_{n-2}),
    & i\leq n-2, i<j\leq n-2.
  \endcases
$$
Let $i\leq n-2,j=n-1$.
In this case we have
$$\align
  \widetilde{D}_{i,n-1}
  &= \frac1{\det D}(-1)^{i+n-1}
     \vmatrix
     -w_1 & 1      & 0        &          &          &   & 0  \\
     1    & \ddots & 1        & 0        &          &   &    \\
     0    & \ddots & -w_{i-1} & 0        &          &   &    \\
          & 0      & 1        & 1        & 0        &   &    \\
          &        & 0        & -w_{i+1} & \ddots   & 0 & 0  \\
          &        & 0        & 1        & -w_{n-2} & 1 & 1  \\
    0     &        &          & 0        & 1        & 0 & -2 \\
    \endvmatrix \\
  &= \frac1{\det D}(-1)^{i+n-1}(-2)a(w_1,\cdots,w_{i-1}).
\endalign$$
Also we have
$$\align
  \widetilde{D}_{in} &= \widetilde{D}_{n-1,i} = \widetilde{D}_{ni} \\
  &= \frac1{\det D}(-1)^{i+n-1}(-2)a(w_1,\cdots,w_{i-1}).
\endalign$$
Now we suppose $i,j\geq n-1$.
In this case, we have
$$\align
  \widetilde{D}_{n-1,n-1}
  &= \widetilde{D}_{n,n} \\
  &= \frac1{\det D}
     \vmatrix
     -w_1 & 1      & 0        & 0  \\
     1    & \ddots & 1        & 0  \\
     0    & 1      & -w_{n-2} & 1  \\
     0    & 0      & 1        & -2 \\
     \endvmatrix
   = \frac1{\det D}a(w_1,\dots,w_{n-2},2).
\endalign$$
$$\align
  \widetilde{D}_{n-1,n}
  &= \widetilde{D}_{n,n-1} \\
  &= \frac1{\det D}(-1)
     \vmatrix
     -w_1 & 1      & 0        & 0 \\
     1    & \ddots & 1        & 0 \\
     0    & 1      & -w_{n-2} & 1 \\
     0    & 0      & 1        & 0 \\
     \endvmatrix
   = \frac1{\det D}a(w_1,\dots,w_{n-3}).
\endalign$$
Note that there is an equation
$$
  d(w_1,\dots,w_{n-2})=-4a(w_1,\dots,w_{n-2},1).
$$
\par
Now we calculate the value $a_i$.
The calculation will proceed as same as of type $A_n$.
But one thing differs, that is $d(w_{n-2})=-4w_2+4$.
For $i\leq n-2$, we have the following.
$$\align
  a_i &= \frac1{\det D}\bigl\{
         d(w_{i+1},\dots,w_{n-2})
           ((-1)^{i+1}-a(w_1,\dots,w_{i-2})-a(w_1,\dots,w_{i-1}))+ \\
      &\hbox to .4in{\hfil}
         a(w_1,\dots,w_{i-1})d(w_{i+1},\dots,w_{n-2})(2-w_i) \\
      &\hbox to .4in{\hfil}
         -a(w_1,\dots,w_{i-1})
           (d(w_{i+1},\dots,w_{n-2})+d(w_{i+2},\dots,w_{n-2}))
         \bigr\} \\
      &= \frac1{\det D}\bigl\{
         (-1)^{i+1}d(w_{i+1},\dots,w_{n-2})
         -a(w_1,\dots,w_{i-2})d(w_{i+1},\dots,w_{n-2}) \\
      &\hbox to .4in{\hfil}
         -a(w_1,\dots,w_{i-1})d(w_{i+1},\dots,w_{n-2})w_i
         -a(w_1,\dots,w_{i-1})d(w_{i+2},\dots,w_{n-2})\bigr\} \\
      &= 1+\frac1{\det D}(-1)^{i+1}d(w_{i+1},\dots,w_{n-2})
       = 1-\frac1{|\det D|}|d(w_{i+1},\dots,w_{n-2})|.
\endalign$$
Since $c_{i,i}=-\widetilde{D}_{ii}$, we have the following.
$$
  a_i+c_{i,i} = 1+\frac1{|\det D|}|d(w_{i+1},\dots,w_{n-2})|
                \bigl(|a(w_1,\dots,w_{i-1})|-1\bigr).
$$
\par
For $i=n-1,n$, we have $a_{n-1}=a_n=a_{n-2}/2$.
Hence we have 
$$
  a_{n-1}=a_n=\frac12+(-1)^{n-1}\frac2{\det D}
         =\frac12-\frac2{|\det D|}.
$$
Therefore we have the following.
$$\align
  a_{n-1}+c_{n-1,n-1} &= a_n+c_{n,n} \\
  &= 1+\frac1{\det D}\{
     2(-1)^{n-1}-\frac12 d(w_1,\dots,w_{n-2})-a(w_1,\dots,w_{n-2},2)\} \\
  &= 1+\frac1{\det D}\{
     2(-1)^{n-1}+2a(w_1,\dots,w_{n-2},1)-a(w_1,\dots,w_{n-2},2)\} \\
  &= 1+\frac1{\det D}\{
     2(-1)^{n-1}-2a(w_1,\dots,w_{n-2})-2a(w_1,\dots,w_{n-3})+ \\
  &\hbox to .4in{\hfil}
     2a(w_1,\dots,w_{n-2})+a(w_1,\dots,w_{n-3})\} \\
  &= 1+\frac1{\det D}\{2(-1)^{n-1}-a(w_1,\dots,w_{n-3})\} \\
  &= 1+\frac1{|\det D|}(|a(w_1,\dots,w_{n-3})|-2) \geq 1.
\endalign$$
Hence we get the following proposition.
\par
\proclaim{Proposition 2}
Let $C$ as above.
If $y$ is of type $D_n$ then $a_i+c_{i,i}\geq 1$ for all $i$.
\endproclaim
\par
\subhead 1.3 The case of type $E_n$ \endsubhead
If the dual graph is in the Appendix then we write \hfil\break
$(m;a,b,c;d,e)$ or may write $(m;a,b,c;d,e;-2)$.
We set indexes $i$ as in the Appendix.
There are only 15 types in all by the classification ([B]).
Hence we can calculate $a_i+c_{i,i}$ directly.
The result is in Appendix.
\par
Hence we get the following proposition.
\par
\proclaim{Proposition 3}
Let $C$ as above.
If $y$ is of type $E_n$ then $a_i+c_i>1$ for all $i$.
\endproclaim
\par
\head 2. Proof of the Theorem \endhead
We recall notation.
Let $(Y,B)$ be a normal surface with a boundary $B$ such that $[B]=0$.
Let $M$ be a Cartier divisor on $Y$.
Let $D$ be an effective $\Q$-divisor on $Y$,
$\Q$-linearly equivalent to the nef big divisor $M-(K_Y+B)$.
Let $y$ be a point of $Y$.
Let $f\:X\to Y$ be the global minimal resolution,
if $y$ is a smooth point, $f$ factors the blowing-up at $y$.
Let $\Delta=f^*K_Y-K_X$.
Let $\Delta_y$ be the components of $\Delta$ supported on $f^{-1}(y)$,
and $\Delta'=\Delta-\Delta_y$.
Let $Z$ be the fundamental cycle of $y$.
Let $\mu=\mu(B,y)$,
where $\mu(B,y)=\max\{\mu\mid\mu(Z-\Delta_y)\leq f^*B\}$.
Let $\delta_y=-(Z-\Delta_y)^2$ if $y$ is quasi-log-terminal, and
$\delta_y=0$ if $y$ is not quasi-log-terminal.
Let $B=B_y+B'$ where $B_y$ consists of the component which contain $y$,
and $B'$ is disjoint from $y$.
Let $f^*B_y=\proper{B_y}+\exc{B_y}$ where $\proper{B}$ is the proper
transformation of $B_y$, and $f(\exc{B_y})=\{y\}$.
Since $M\sim K_Y+D+B$ is Cartier,
$f^*(K_Y+D+B)=K_X+f^*(D+B)+\Delta$ is integral.
Hence $D+B$ is also integral coefficients, thus $\rup{D}=D+B$.
\par
First we treat the case which $y$ is not quasi-log-terminal.
In this case, we assume that $D^2>0$ and
$DC\geq 0$ for all curve $C\ni y$, because we defined $\delta_y=0$.
Since $D$ is nef and big, this always hold.
Hence it is sufficient to prove that $y$ is not a base point of $|M|$.
\par
Since $f^*D$ is nef and big, we have $H^1(X,K_X+\rup{f^*D})=0$
by Sakai's lemma [S1] where $\rup{}$ means the ``round-up''.
On the other hand,
$$\align
  K_X+\rup{f^*D}
  &= \rup{K_X+f^*D}=\rup{f^*(K_Y+B)-\Delta-f^*B+f^*D} \\
  &= \rup{f^*(K_Y+\rup{D})-(\Delta_y+\exc{B_y})-\Delta'-f^*B'} \\
  &= f^*(K_Y+\rup{D})-\rdown{\Delta_y+\exc{B_y}}-\rdown{\Delta'+f^*B'}.
\endalign$$
Since $y$ is not quasi-log-terminal,
the round-down $\rdown{\Delta_y+\exc{B_y}}\ne 0$.
Let $S=\rdown{\Delta_y+\exc{B_y}}$ and $T=\rdown{\Delta'+f^*B'}$.
Hence we have a surjection
$$
  H^0(X,f^*(K_Y+\rup{D})-T)\to H^0(S, (f^*(K_Y+\rup{D})-T)|_S)\to 0.
$$
Since $f(S)=\{y\}$, $\O_S((f^*(K_Y+\rup{D})-T)|_S)$ is trivial.
Hence we get a section of $\O_X(f^*(K_Y+\rup{D})-T)$ which does not vanish
at any point of $\Supp S$.
Therefore we get a global section of $\O_X(f^*(K_Y+\rup{D}))$
which does not vanish on $f^{-1}(y)$ by multiplying
the global section of $\O_X(T)$.
\par
From now we assume that $y$ is quasi-log-terminal.
Let $D=D_y+D'$ where all irreducible divisors in $D_y$ contain $y$
and $D'$ is the others.
Let $f^*D_y=\proper{D_y}+\exc{D_y}=\sum d_i D_i+\sum d_j'\Delta_j$.
Let $f^*B_y=\sum b_i D_i+\sum b_j'\Delta_j$.
Let $\Delta_y=\sum a_j\Delta_j$.
Since $D^2> -((1-\mu)(Z-\Delta_y))^2$, $f^*D-(1-\mu)(Z-\Delta_y)$ is big.
Hence if necessary, we change $D$ by $\Q$-linearly equivalence and
we assume $D_y\ne 0$.
We define the number 
$$
  c=\min\left\{\frac{1-b_i}{d_i},\frac{1-a_j-b_j'}{d_j'}\right\}
$$
Since $d_i,d_j'>0$, $0\leq b_i<1$ and $a_j+b_j'<1$ we have $c>0$.
\par
\proclaim{Lemma 3}
$c\leq 1/2$.
\endproclaim
\par
\demo{Proof}
Since $D^2>-((1+\varepsilon)(1-\mu)(Z-\Delta_y))^2$
for sufficiently small $\varepsilon$,
we have $D-(1+\varepsilon)(1-\mu)(Z-\Delta_y)>0$.
Hence we have $d_j'-(1-\mu)(z_j-a_j)>0$ where $Z=\sum z_j\Delta_j$.
Thus $d_j'+a_j-z_j+\mu(z_j-a_j)>0$.
By the definition of $\mu$, we have $b_j'\geq \mu(z_j-a_j)$.
Hence we have $d_j'+a_j-z_j+b_j'>0$.
Therefore $d_j'+a_j+b_j'>z_j\geq 1$ because $z_j$ is positive integer.
On the other hand, Since $K_Y+\rup{D}=K_Y+D+B$ is Cartier,
$f^*D+f^*B+\Delta=f^*(K_Y+D+B)-K_X$ is integral divisor.
Hence $d_j'+a_j+b_j'$ is also integral.
Since $d_j'+a_j+b_j'\geq 2$, we have $2(1-a_j/2-b_j'/2)\leq d_j'$.
Thus we have
$$
  c\leq\frac{1-a_j-b_j'}{d_j'}\leq\frac12-\frac{a_j+b_j'}{2d_j'}
   \leq\frac12. \qed
$$
\enddemo
\par
Let
$$\align
  R &= f^*D-cf^*D=f^*\rup{D}+\Delta-f^*B-cf^*D-\Delta \\
    &= f^*\rup{D}+\Delta
       -\left(\sum (b_i+cd_i)D_i+\sum (b_j'+cd_j'+a_j)\Delta_j\right)
       -(f^*B'+f^*D'+\Delta').
\endalign$$
Since $f^*\rup{D}+\Delta=f^*(K_Y+\rup{D})-K_X$ is Cartier,
$$\align
  \rup{R}
    &= f^*\rup{D}+\Delta
       -\rdown{\sum (b_i+cd_i)D_i+\sum (b_j'+cd_j'+a_j)\Delta_j}
       -\rdown{f^*B'+f^*D'+\Delta'} \\
    &= f^*\rup{D}+\Delta-(D_1+\dots +D_s)-E'-Q
\endalign$$
where $\rdown{\sum (b_i+cd_i)D_i}=D_1+\cdots +D_s$,
$\rdown{\sum (b_j'+cd_j'+a_j)\Delta_j}=E'$, and
\hfill\break  
$\rdown{f^*B'+f^*D'+\Delta'}=Q$.
\par
Suppose $E'\ne 0$.
Since $R=(1-c)f^*D$ is nef big and $D_1,\dots,D_s$ are not appear
in the fractional part of $R$,
we have $H^1(X,K_X+\rup{R}+D_1+\cdots+D_s)=0$ by
Kawamata-Viehweg vanishing theorem.
Hence we have the surjection
$$
  H^0(X,f^*(K_Y+\rup{D})-Q)\to H^0(E',(f^*(K_Y+\rup{D})-Q)|_{E'})\to 0.
$$
Since $f(E')=\{y\}$ and $Q$ is disjoint from $y$,
we get a global section of $K_Y+\rup{D}$ which is not vanish at $y$.
\par
Now we assume $E'=0$.
Note that the minimality of $c$ induce that $b_i+cd_i\leq 1$.
Thus the minimal value $c$ is obtained by $(1-b_i)/d_i$,
hence $s\geq 1$.
By Kawamata-Viehweg vanishing theorem, we have
$$
  H^0(X, K_X+\rup{R}+D_2+\cdots+D_s)= H^0(K_X+f^*\rup{D}+\Delta-D_1-Q)=0.
$$
Thus there is the surjection
$$
  H^0(X,K_X+f^*\rup{D}+\Delta-Q)\to
  H^0(D_1,K_{D_1}+(f^*\rup{D}+\Delta-D_1-Q)|_{D_1})\to 0.
$$
So it is sufficient to prove $\rup{R}D_1>1$.
Indeed, if $\rup{R}D_1>1$ then $\rup{R}D_1\geq 2$
because $\rup{R}$ is Cartier.
Thus $(f^*\rup{D}+\Delta-D_1-Q)D_1\geq\rup{R}D_1\geq 2$.
Hence we get the nowhere vanishing global section in
$H^0(D_1,K_{D_1}+(f^*\rup{D}+\Delta-D_1-Q)|_{D_1})$.
Hence we get the global section of $\O(f^*(K_Y+\rup{D}))$
which is not vanish on $D_1$, especially not vanish at $y$.
\par
Hence we prove $\rup{R}D_1>1$.
Remark that
$$
  \rup{R}=R+\left\{\sum_{i\ne 1,\dots s}(b_i+cd_i)D_i+
                   \sum (b_j'+cd_j'+a_j)\Delta_j+\text{(others)}\right\}
$$
where $\{\text{(divisor)}\}$ is the fractional part of (divisor).
Since $E'=0$, 
$\{\sum (b_j'+cd_j'+a_j)\Delta_j\}=\sum (b_j'+cd_j'+a_j)\Delta_j$.
Therefore
$$
  \rup{R}D_1\geq (1-c)f^*DD_1+\sum (b_j'+cd_j'+a_j)\Delta_jD_1.
$$
Thus it is sufficient to prove
$$
  (1-c)(1-\mu)\delta_y/2+\sum (b_j'+cd_j'+a_j)\Delta_jD_1>1. \tag 1
$$
\par
If $D_1$ meets more components of $\Delta_y$ then
$\sum (b_j'+cd_j'+a_j)\Delta_jD_1$ is larger.
Thus we assume $\Delta_iD_1=1$ for an $i$ and $\Delta_iD_1=0$ for $j\ne i$.
Let $f^*f_*D_1=D_1+\sum c_{i,j}\Delta_j$.
Since $f^*B_y+cf^*D_y=D_1+\text{(others)}$, we have $c_{i,j}\leq b_j'+cd_j'$.
But by Proposition 1 through 3, $c_{i,i}+a_i\geq 1$ if $y$ is of type $D_n$
or of type $E_n$, and if $y$ is of type $A_n$ and $i\ne 1,n$.
They contradict to the assumption $E'=0$.
Thus we assume $y$ is of type $A_n$ and $j=1,n$ or $y$ is a smooth point.
\par
Thus if $y$ is of type $D_n$ or $E_n$ and $D^2>(1-\mu)^2\delta_y$
then $y$ is not a base point of $|M|$.
\par
First we assume that $y$ is a smooth point or of type $A_1$.
In this case, $b_1'=\mu(1-a_1)$.
Let $p=\mu(D,y)$, hence $d_1'=p(1-a_1)$.
Let $w_1=-\Delta_1^2$.
Then we have $a_1=1-2/w_1$ and $\delta_y=4/w_1$.
Hence $\delta_y=2(1-a_1)$.
Since $f^*D-(1-\mu)(Z-\Delta_y)$ is big, we have $p+\mu>1$
as in the proof of lemma 3.
Hence we have
$$\align
  \rup{R}D_1
    &\geq (1-c)(1-\mu)\frac{\delta_y}2+(b_1'+cd_1'+a_1) \\
    &= (1-\mu)\frac{\delta_y}2+b_1'+a_1+(d_1'-(1-\mu)\frac{\delta_y}2)c \\
    &= (1-\mu)(1-a_1)+\mu(1-a_1)+a_1+(p(1-a_1)-(1-\mu)(1-a_1))c \\
    &= 1+(p+\mu-1)(1-a_1)c >1. \\
\endalign$$
\par
Next we assume that $y$ is of type $A_n$ and $n\geq 2$.
We may change the indexes of $\Delta_j$, we assume $a_1\leq a_n$.
We also assume that $D_1$ meets $\Delta_i$ where $i=1$ or $n$.
\par
Case 1: We assume that $i=n$ or $a_1=a_n$.
\par
It is  sufficient to prove $(1-c)(1-\mu)\delta_y/2+b_i'+cd_i'+a_i>1$
by the inequality of (1).
Since the definition of $\mu$, we have $\mu(1-a_i)\leq b_i'$.
Since $f^*D-(1-\mu)(Z-\Delta_y)$ is big,
we have $d_i'> (1-\mu)(1-a_i)$.
Since $\delta_y=2-a_1-a_n$, we have
$$\align
  (1 &- c)(1-\mu)\delta_y/2+b_i'+cd_i'+a_i \\
    &> (1-\mu)\left(1-\frac{a_1+a_n}2\right)+\mu(1-a_i)+a_i+ \\
    &\quad \left((1-\mu)(1-a_i)
           -(1-\mu)\left(1-\frac{a_1+a_n}2\right)\right)c \\
    &= 1-\frac{a_1+a_n}2+a_i+\mu\left(1-a_i-1+\frac{a_1+a_n}2\right)+
       c(1-\mu)\left(1-a_i-1+\frac{a_1+a_n}2\right) \\
    &= 1+\frac{a_n-a_1}2-\mu\frac{a_n-a_1}2-c(1-\mu)\frac{a_n-a_1}2
       \quad \text{(let $i=n$)} \\
    &= 1+(1-c)(1-\mu)\frac{a_n-a_1}2 \geq 1.
\endalign$$
\par
Case 2: We assume that $i=1$ and $a_1<a_n$.
\par
Let $B_y=bD_1+B''$ and $\mathop{ord}_{\Delta_i} f^*B''=b_i''$.
Since
$$
  \sum (b_i'+cd_i'+a_i)\Delta_i=\exc{(f^*B_y+cf^*D_y+\Delta_y)}
  \geq \exc{(f^*f_*D_1+f^*B''+\Delta_y)},
$$
we have $b_i'+cd_i+a_i\geq a_i+c_i+b_i''$.
Since $a_1+c_{1,1}=1-1/|\det A|$,
if $b_1''\geq 1/|\det A|$ then $\rup{R}D_1>1$ because $\delta_y>0$.
\par
Suppose that $b_1''<1/|\det A|$.
Let $B''=\sum b_k D_k$ and
$f^*D_k=\proper{D_k}+\sum (c_{k_1,j}+\cdots+c_{k_{r_k},j})\Delta_j$
if $\proper{D_k}$ meets $\Delta_{k_1},\dots,\Delta_{k_{r_k}}$
counted with intersection multiplicity.
Since 
$$
  \frac1{|\det A|}>b_1''=\sum_{k,i} b_kc_{k_i,1}
  \geq \sum_k b_kc_{n,1}r_k=\frac{\sum b_kr_k}{|\det A|},
$$
we have $\beta=\sum b_k r_k<1$.
By the definition of $\mu$ we have
$$\align
  \mu &\leq \frac{b_n''}{1-a_n}
      = \frac{bc_{1,n}+\sum b_kc_{k_i,n}}{1-a_n}
      \leq \frac{bc_{1,n}+\sum b_kr_kc_{n,n}}{1-a_n} \\
      &= \frac{bc_{1,n}+\beta c_{n,n}}{1-a_n}.
\endalign$$
Hence we have
$$\align
  (1-c) & (1-\mu)\frac{\delta_y}2 \\
  &\geq (1-c)\left(1-\frac{bc_{1,n}+\beta c_{n,n}}{1-a_n}\right)
        \frac{\delta_y}2 \\
  &= (1-c)\frac{\delta_y}{2(1-a_n)}
       ((1-a_n)|\det A|-b-\beta a(w_1,\dots,w_{n-1}))\frac1{|\det A|} \\
  &= (1-c)\frac{\delta_y}{2(1-a_n)}
       ((1-b)+(1-\beta)a(w_1,\dots,w_{n-1}))\frac1{|\det A|}.
\endalign$$
If $b\leq \beta$ then
$$\align
  (1-c)(1-\mu)\frac{\delta_y}2
  &\geq (1-c)\frac{\delta_y}{2(1-a_n)}(1-\beta)
        \frac{1+a(w_1,\dots,w_{n-1})}{|\det A|} \\
  &= (1-c)(1-\beta)\frac{\delta_y}{2(1-a_n)}(1-a_n) \\
  &= (1-c)(1-\beta)\frac{\delta_y}2.
\endalign$$
Since $c\leq 1/2$ by lemma 3, we have
$$\align
  (1-c)& (1-\mu)\frac{\delta_y}2+\frac\beta{|\det A|} \\
  &\geq (1-\beta)\frac{\delta_y}4+\frac\beta{|\det A|} \\
  &= \frac1{|\det A|}\left(\frac14(2+a(w_1,\dots,w_{n-1})
     +a(w_2,\dots,w_n))(1-\beta)+\beta\right) \\
  &\geq \frac1{|\det A|}\left(\frac32-\frac12\beta\right) \\
  &> \frac1{|\det A|}.
\endalign$$
Hence we have
$$\align
  (1-c)& (1-\mu)\frac{\delta_y}2+b_1'+cd_1+a_1 \\
  &\geq (1-c)(1-\mu)\frac{\delta_y}2+a_1+c_{1,1}+\sum b_k c_{k_i,1} \\
  &\geq (1-c)(1-\mu)\frac{\delta_y}2+1-\frac1{|\det A|}+\frac\beta{|\det A|} \\
  &> 1.
\endalign$$
\par
Finally we assume that $\beta< b$.
In this case we have
$$\align
  (1-c) & (1-\mu)\frac{\delta_y}2 \\
  &\geq \frac{\delta_y}{4(1-a_n)}
       ((1-b)+(1-\beta)a(w_1,\dots,w_{n-1}))\frac1{|\det A|} \\
  &= \frac{\delta_y}{4(1-a_n)}(1-\beta)(1+a(w_1,\dots,w_{n-1}))\frac1{|\det A|}
     -\frac{\delta_y}{4(1-a_n)}\frac1{|\det A|}(b-\beta) \\
  &> (1-\beta)\frac{\delta_y}4
     -\frac{\delta_y}4\frac1{1+a(w_1,\dots,w_{n-1})}(1-\beta) \\
  &= (1-\beta)\frac{\delta_y}4\left(1-\frac1{1+a(w_1,\dots,w_{n-1})}\right) \\
  &\geq (1-\beta)\frac{\delta_y}6.
\endalign$$
Hence we have
$$\align
  (1-c)(1-\mu)\frac{\delta_y}2 &+ \frac\beta{|\det A|} \\
  &> (1-\beta)\frac{\delta_y}6+\frac\beta{|\det A|} \\
  &= \frac1{|\det A|}\left(\frac16(2+a(w_1,\dots,w_{n-1})
     +a(w_2,\dots,w_n))(1-\beta)+\beta\right) \\
  &\geq \frac1{|\det A|}.
\endalign$$
Therefore we have
$$
  (1-c)(1-\mu)\frac{\delta_y}2+b_1'+cd_1+a_1 >1.
$$
Thus we have $\rup{R}D_1>1$, we complete the proof. \qed
\par
\newpage
\head Appendix \endhead
$$\vbox{\offinterlineskip\let\sc\scriptstyle\let\frac\tfrac
  \halign{\vrule#&\strut\ \hfil#\hfil &\vrule#&\ $\sc#$\hfil\ &\vrule#&
          \ $\sc#$\hfil\ &\vrule#&
          \ $\sc#$\hfil\vbox{\hrule width0pt height 11pt}\ &\vrule# \cr
    \noalign{\hrule} 
    & type && \text{\tenrm dual graph} && a_i+c_i && x= & \cr
    height2pt & \omit && \omit && \omit && \omit &\cr
    \noalign{\hrule} 
    & 1 && (m;2,2;2,2) &&
           \left(1+\frac5x;1+\frac1{3x},\frac43+\frac2x;\frac43+\frac2x,
           1+\frac1{3x};1+\frac1x\right) && 6m-11 & \cr
    height2pt & \omit && \omit && \omit && \omit & \cr
    \noalign{\hrule}
    & 2 && (m;2,2;3) && 
           \left(1+\frac5{3x};1+\frac1{9x},\frac43+\frac2{3x};
           1+\frac1{9x};1+\frac1{3x}\right) && 2m-3 & \cr
    height2pt & \omit && \omit && \omit && \omit & \cr
    \noalign{\hrule}
    & 3 && (m;3;3) && 
           \left(1+\frac5x;1+\frac1{3x};1+\frac1{3x};1+\frac1x
           \right) && 6m-7 & \cr
    height2pt & \omit && \omit && \omit && \omit & \cr
    \noalign{\hrule}
    & 4 && (m;2,2,2;2,2) && 
           \left(1+\frac{11}x;1+\frac1{2x},\frac32+\frac5{2x},
           \frac32+\frac6x;\frac43+\frac{14}{3x},1+\frac1x;
           1+\frac5{2x}\right) && 12m-23 & \cr
    height2pt & \omit && \omit && \omit && \omit & \cr
    \noalign{\hrule}
    & 5 && (m;2,2,2;3) && 
           \left(1+\frac{11}x;1+\frac1{2x},\frac32+\frac5{2x},
           \frac32+\frac6x;1+\frac1x;1+\frac5{2x}\right) && 12m-19 & \cr
    height2pt & \omit && \omit && \omit && \omit & \cr
    \noalign{\hrule}
    & 6 && (m;4;2,2) && 
           \left(1+\frac{11}x;1+\frac1{2x};\frac43+\frac14{3x},
           1+\frac1x;1+\frac5{2x}\right) && 12m-17 & \cr
    height2pt & \omit && \omit && \omit && \omit & \cr
    \noalign{\hrule}
    & 7 && (m;4;3) && 
           \left(1+\frac{11}x;1+\frac1{2x};1+\frac1x;
           1+\frac5{2x}\right)&& 12m-13 & \cr
    height2pt & \omit && \omit && \omit && \omit & \cr
    \noalign{\hrule}
    & 8 && (m;2,2,2,2;2,2) &&
           \left(1+\frac{29}x;1+\frac1x,\frac85+\frac{22}{5x},
           \frac95+\frac{51}{5x},\frac85+\frac{92}{5x};
           \frac43+\frac{38}{3x},1+\frac3x;
           1+\frac7x\right) && 30m-59 & \cr
    height2pt & \omit && \omit && \omit && \omit & \cr
    \noalign{\hrule}
    & 9 && (m;2,2,2,2;3) && 
           \left(1+\frac{29}x;1+\frac1x,\frac85+\frac{22}{5x},
           \frac95+\frac{51}{5x},\frac85+\frac{92}{5x};
           1+\frac3x;1+\frac7x\right) && 30m-49 & \cr
    height2pt & \omit && \omit && \omit && \omit & \cr
    \noalign{\hrule}
    & 10 && (m;2,3;2,2) && 
            \left(1+\frac{29}x;1+\frac1x,\frac65+\frac{22}{5x};
            \frac43+\frac{38}{3x},1+\frac3x;
            1+\frac7x\right) && 30m-47 & \cr
    height2pt & \omit && \omit && \omit && \omit & \cr
    \noalign{\hrule}
    & 11 && (m;2,3;3) && 
            \left(1+\frac{29}x;1+\frac1x,\frac65+\frac{22}{5x};
            1+\frac3x;1+\frac7x\right) && 30m-37 & \cr
    height2pt & \omit && \omit && \omit && \omit & \cr
    \noalign{\hrule}
    & 12 && (m;3,2;2,2) && 
            \left(1+\frac{29}x;1+\frac1x,\frac75+\frac{51}{5x};
            \frac43+\frac{38}{3x},1+\frac3x;
            1+\frac7x\right) && 30m-53 & \cr
    height2pt & \omit && \omit && \omit && \omit & \cr
    \noalign{\hrule}
    & 13 && (m;3,2;3) && 
            \left(1+\frac{29}x;1+\frac1x,\frac75+\frac{51}{5x};
            1+\frac3x;1+\frac7x\right) && 30m-43 & \cr
    height2pt & \omit && \omit && \omit && \omit & \cr
    \noalign{\hrule}
    & 14 && (m;5;2,2) && 
            \left(1+\frac{29}x;1+\frac1x;\frac43+\frac{38}{3x},
            1+\frac3x;1+\frac7x\right) && 30m-41 & \cr
    height2pt & \omit && \omit && \omit && \omit & \cr
    \noalign{\hrule}
    & 15 && (m;5;3) && 
            \left(1+\frac{29}x;1+\frac1x;1+\frac3x;
            1+\frac7x\right) && 30m-31 & \cr
    height2pt & \omit && \omit && \omit && \omit & \cr
    \noalign{\hrule}
    }}
$$
\par
\newpage
\topinsert
\vskip 8truein
\includegraphics{DualGraph.ps}
\endinsert
\vbox to 0pt{}
\newpage
\enddocument